\documentclass[11pt]{amsart}

\usepackage{amsfonts, amsmath, amssymb}
\usepackage{graphicx, graphics}

\newtheorem{prop}{Proposition}[section]
\newtheorem{theo}[prop]{Theorem}
\newtheorem{coro}[prop]{Corollary}

\newtheorem{conj}[prop]{Conjecture}

\theoremstyle{definition}
\newtheorem{dfn}[prop]{Definition}

\theoremstyle{remark}
\newtheorem{rem}[prop]{Remark}

\def\<{\langle}
\def\>{\rangle}

\def\G{{\Gamma}}

\def\O{{\mathcal O}}
\def\P{{\mathbb P}}

\def\Z{{\mathbb Z}}

\def\A{{\mathbb A}}

\newcommand{\OO}{\mathcal{O}}
\newcommand{\KK}{\mathcal{K}}
\newcommand{\ZZ}{\mathbb{Z}}
\newcommand{\PP}{\mathbb{P}}

\newcommand{\comment}[1]{}

\newcommand{\pp}{$/\!/$ }

\theoremstyle{definition}

\DeclareMathOperator{\Pic}{Pic}

\DeclareMathOperator{\Proj}{Proj}
\DeclareMathOperator{\Spec}{Spec}
\DeclareMathOperator{\Cox}{Cox}

\newsavebox{\fourpts}
\savebox{\fourpts}(1.5,1.5)[bl]
{\small\tt
\multiput(0,0)(1.5,0){2}{\circle*{.1}}
\multiput(0,1.5)(1.5,0){2}{\circle*{.1}}
\put(0.1,.1){1}
\put(1.6,.1){2}
\put(0.1,1.6){3}
\put(1.6,1.6){4}
}

\begin{document}

\title{The Cox ring of a Del Pezzo surface}

\author{Victor~V.~Batyrev}
\address{Mathematisches Institut Universit\"at T\"ubingen, 
Auf der Morgenstelle~10, T\"ubingen D-72076, Germany} 
\email{victor.batyrev@uni-tuebingen.de} 

\author{Oleg~N.~Popov}
\address{Department of Algebra, Faculty of Mathematics, Moscow State 
University, Moscow 117234, Russia}  
\email{popov@mccme.ru}


\keywords{Del Pezzo surfaces, torsors, homogeneous spaces, algebraic groups} 

\begin{abstract} 
Let $X_r$ be  a smooth Del Pezzo surface obtained from $\P^2$ 
by blowing up  $r \leq 8$ points in general position. 
It is well known that for $r \in \{ 3,4,5,6,7,8 \}$ the Picard group 
$\Pic(X_r)$ contains a canonical  root system 
$R_r \in \{ A_2 \times A_1, A_4, D_5, E_6, E_7, E_8 \}$. 
We prove some general properties of the  Cox ring of 
$X_r$ ($r \geq 4$) and  show its similarity to the homogeneous 
coordinate ring of the orbit of the highest weight vector 
in some irreducible representation of the algebraic group $G$ associated 
with the root system $R_r$. 
\end{abstract}

\maketitle

\section{Introduction}

Let $X$ be a  projective algebraic variety over a field 
$\Bbbk$. Assume  that the Picard group 
$\Pic(X)$ is a finitely generated abelian group. Consider 
the vector space 
\[ \G(X):= \bigoplus_{ [D] \in \Pic(X)} H^0(X, \O(D)). \]
One wants to make it an $\Bbbk$-algebra which is  graded by the monoid of 
effective classes in  $\Pic(X)$ such that the algebra structure will 
be compatible with the natural bilinear map 
\[ b_{D_1,D_2}\;: \;  H^0(X, \O(D_1)) \times  
H^0(X, \O(D_2)) \to  H^0(X, \O(D_1 + D_2)). \]
However, there exist some  problems in the realization of this idea. 
First of all 
 there is no any natural isomorphism between 
$H^0(X, \O(D))$ and  $H^0(X, \O(D'))$ if $[D] = [D']$. 
There exists only a canonical bijection   
between the linear systems $|D| \cong  |D'|$ (where $|D|$ is the  
projectivization  of the $\Bbbk$-vector space  $H^0(X, \O(D))$). 
As a consequence, the bilinear map 
$b_{D_1,D_2}$  depends not only on the classes  $[D_1], [D_2], 
[D_1 + D_2] \in \Pic(X)$, but also on their particular representatives.  
One can easily see that only the  morphism   
\[  s_{[D_1],[D_2]}\;: \; |D_1| \times  
|D_2| \to  |D_1 + D_2|  \]
of the product of two projective spaces  $|D_1| \times  
|D_2|$ to another projective space $|D_1 + D_2|$ is well-defined. 
For this reason, it is much more natural to consider 
the graded set of projective spaces  
\[ |\G(X)|:= \bigsqcup_{ [D] \in \Pic(X)} |D| \]
together with all possible morphisms  $s_{[D_1],[D_2]}$ any two effective 
classes $[D_1],[D_2] \in 
\Pic(X)$. 

Inspired by the paper of  Cox  on the homogeneous ring of a toric 
variety \cite{Cox}, Hu and Keel \cite{HK} suggested a definition 
of a {\em Cox ring} 
\[\Cox(X) = R(X, L_1, \ldots, L_r) := \bigoplus_{(m_1, \ldots, m_r) \in \Z^r} 
H^0(X, \O(m_1L_1 + \cdots + m_rL_r)) \]
which uses a choice of some $\Z$-basis $L_1, \ldots, L_r$ in  $\Pic(X)$ 
(e.g. if $\Pic(X) \cong \Z^r$ is a free abelian group). 
Using  such a $\Z$-basis, one 
obtains a particular representative for each class in  $\Pic(X)$ together   
with a well-defined multiplication
so  $R(X, L_1, \ldots, L_r)$ becomes a well-defined $\Bbbk$-algebra. 
If $L_1', \ldots, L_r'$ is another $\Z$-basis of $\Pic(X)$, then 
the corresponding Cox algebra  $R(X, L_1', \ldots, L_r')$ is 
isomorphic to $R(X, L_1, \ldots, L_r)$. 
Unfortunately, we can not expect to choose a  $\Z$-basis
of $\Pic(X)$ in a natural canonical way. More often one 
can  choose in a  natural  way some 
effective divisors $D_1, \ldots, D_n$ on $X$ such that $\Pic(X)$ 
is generated by 
$[D_1], \ldots, [D_n]$. If we set 
$$U:= X \setminus (D_1 \cup \cdots \cup D_n)$$
 and assume that 
$X$ is smooth, then $Pic(U) = 0$ and we obtain the 
exact sequence 
\[ 1 \to \Bbbk^* \to  \Bbbk[U]^* \to \bigoplus_{i =1}^n \Z[D_i] \to 
\Pic(X) \to 0. \]

Choosing a $\Bbbk$-rational point $p$ in $U$, we can split the monomorphism
$\Bbbk^* \to  \Bbbk[U]^*$, so that one has an isomorphism 
\[ \Bbbk[U]^* \cong \Bbbk^* \oplus G, \]
where $G \subset \Bbbk[U]^*$ is a free abelian group of rank $n-r$.   
The choice of a $\Bbbk$-rational point $p \in U$ allows to give 
another approach to the graded space $\G(X)$ and to the Cox algebra:

\begin{dfn} Let $X, U, p,  D_1, \ldots, D_n$ be as above. 
We consider  the graded $\Bbbk$-algebra  
\[ \G(X, U, p) := \bigoplus_{(m_1, \ldots, m_n) \in \Z^n} 
H^0(X, \O(m_1D_1 + \cdots + m_nD_n)) \]
and define 
\[\Cox(X,U,p):= \G(X, U,p)_G \]
as the quotient of the  $\G(X, U,p)$ modulo the ideal 
generated by 
\[ \{ x - gx \; | \; x \in  \G(X, U,p), g \in G \}. \] 
Since $Pic(X) \cong \Z^n/G$, we obtain a natural 
$\Pic(X)$-grading on  $\Cox(X,U,p)$. 
\label{cox-ring}
\end{dfn} 

We expect that the algebra $\Cox(X,U,p)$ can de applied to 
some arithmetic questions on $\Bbbk$-rational points in  $U \subset X$.

\begin{rem} The above definition  of the ring  $\Cox(X,U,p)$ depends on 
the choice of an  open subset $U \subset X$ and a $\Bbbk$-rational point
$p \in U$.  A similar idea was used by Colliot-Th{\'e}l{\`e}ne and Sansuc in 
\cite{CS} for 
constructing  universal torsors and deriving  explicit equations for them. 
The lack of a canonical construction is precisely what makes 
descending the universal torsor an interesting problem. 
Some   applications of the universal torsor for Del Pezzo surfaces 
of degree $5$  was  considered by Skorobogatov 
in \cite{Sk1} (see also \cite{Sk2}). Recently, Hassett and  Tschinkel have 
investigated the Cox rings and the universal torsors  for some 
interesting  special cubic 
surfaces \cite{HT}. 
\end{rem}

\begin{rem} 
If $X$ is a smooth projective toric variety and $U \subset X$ is the open 
dense torus orbit, then the choice of a point $p \in U$ defines an 
isomorphism of $U$ with the algebraic torus $T$, so that the subgroup 
$G \subset  \Bbbk[U]^*$ can be identified with the character group 
of $T$. In this way, one can show that  $\Cox(X,U,p)$ is isomorphic to 
a polynomial ring in $n$ variables $(n $ is the number of irreducible 
components of $X \setminus U$, cf. \cite{Cox}). 
\label{cox-tor}
\end{rem} 

\begin{rem} 
The field of fractions of the ring   $\Cox(X,U,p)$ is  
a pure  transcendental extension of degree $r$ of  the field 
of rational functions on $X$. 
Therefore, $\dim\, \Spec \,  \G(X, U, p) = \dim X + r$, if $\G(X, U, p)$ 
is a finitely generated $\Bbbk$-algebra. 
\label{dimension}
\end{rem} 

\bigskip

Let $X_r$ be a smooth Del Pezzo surface obtained from $\P^2$ 
by blow-up of $r \leq 8$ points in general position. 
It is well known that for $r \in \{ 3,4,5,6,7,8 \}$ the Picard group 
$\Pic(X_r)$ contains a canonical  root system 
$R_r \in \{ A_2 \times A_1, A_4, D_5, E_6, E_7, E_8 \}$. 
Moreover, the natural embedding $\Pic(X_{r-1}) \hookrightarrow \Pic(X_{r})$ 
 induces the inclusion of root systems $R_{r-1}  \hookrightarrow R_r$. 
If $G(R_r)$ is a connected algebraic group corresponding to the root 
system $R_r$, then the embedding  $R_{r-1}  \hookrightarrow R_r$ 
defines a maximal parabolic subgroup $P(R_{r-1}) \subset G(R_r)$ \cite{Hm}. 
We expect that for $r \geq 4$ there should be some relation 
between a Del Pezzo surface 
$X_r$ and the GIT-quotient of the homogeneous space  $ G(R_r)/P(R_{r-1})$ 
modulo the action of a  maximal torus $T_r$ of   $G(R_r)$.

Our starting 
observation  is the well-known isomorphism $X_4 \cong G(3,5)//T_4$ 
which follows
from an isomorphism between the homogeneous coordinate ring of the 
Grassmaniann $G(3,5) = G(A_4)/P(A_2 \times A_1) \subset \P^9$ and  
the Cox ring of $X_4$ (see \ref{G35}). Another
proof of this fact follows form the identification of $X_4$ with 
the moduli space  $\overline{M_{0,5}}$ of stable rational curves with 
$5$ marked points \cite{Kapranov}.

In this paper, we start an investigation of  the Cox ring of 
Del Pezzo surfaces $X_r$ 
$(r \geq 4)$. It is natural to choose the  classes of 
all exceptional curves 
$E_1, \ldots, E_{N_r} \subset X_r$ as a generating
set  for the Picard group $\Pic(X_r)$.
There is a natural $\Z_{\geq 0}$-grading on  $\Pic(X_r)$ defined by the 
intersection with the anticanonical divisor $-K$. 

We prove some general properties of the  Cox rings of a Del Pezzo 
surface $X_r$ ($r \geq 4$) and  show their similarity to the homogeneous 
coordinate ring of $G(R_r)/P(R_{r-1})$. We remark that  the homogeneous 
space $G(R_r)/P(R_{r-1})$ 
can be interpreted as the  orbit of the highest weight vector in some 
natural irreducible representation of $G(R_r)$.    

\begin{rem} 
Some other connections between  Del Pezzo surfaces and the corresponding
algebraic groups were considered also by Friedman and Morgan in~\cite{F-M}. 
A similar topic was considered by Leung in~\cite{Le}.
\end{rem} 

In this paper, we show that the Cox ring of a Del Pezzo surface $X_r$ 
is generated by elements of degree $1$.  This implies that 
the homogeneous coordinate ring of $G(R_r)/P(R_{r-1})$ is naturally 
graded by  the monoid of effective divisor classes on the surface $X_r$ 
(the same monoid defines the multigrading of the Cox ring of $X_r$).  
Moreover, we obtain some results of the quadratic
relations between the generators of the Cox ring of $X_r$. 
\bigskip

The authors would like to thank Yu. Tschinkel, A. Skorobogatov,  
E.~S.~Golod, S.~M.~Lvovski and E.~B.~Vinberg for useful discussions and
encouraging remarks.

\section{Del Pezzo Surfaces}

Let us summarize briefly some well-known classical results on Del Pezzo 
surfaces which can be found  in 
\cite{Ma,Dem,Na}.

One says  that $r$ $(r \leq 8)$ points $p_1, \ldots, p_r $ in $\P^2$  
are in {\em general position} if there are no 3 points on
a line, no 6 points on a conic ($r\geq 6$) and a cubic having seven 
points and one of
them double does not have the eighth one ($r=8$). 

Denote by $X_r$ $(r \geq 3)$ 
the Del Pezzo surfaces obtained from $\P^2$ by blowing up 
of $r$ points  $p_1, \ldots, p_r$ in general position. 
If  $\pi: X_r\to \P^2$ the corresponding projective morphism, then 
the Picard group $Pic(X_r) \cong \Z^{r+1}$ contains a $\Z$-basis 
$l_i, (0\leq i\leq r)$, $l_0=[\pi^*{\mathcal O}(1)]$ and 
$l_i:=[\pi^{-1}(p_i)],\,  i=1,\dots,r$. The
intersection  form $(*,*)$ on  $Pic(X_r)$ is determined in the chosen basis 
by the diagonal matrix:
$(l_0,l_0)=1, (l_i,l_i)=-1$ for $i\geq 1,$ 
$(l_i,l_j)=0$ for $i\ne j$. The anticanonical class of $X_r$ equals
$-K=3l_0-l_1-\dots-l_r$. The number  $d:= (K,K) = 9 -r$ is called the 
{\em degree} of $X_r$. 
The anticanonical system $|-K|$  of a Del Pezzo surface $X_r$ 
 is very ample if $r\leq 6$,
 it determines a two-fold covering of  $\PP^2$ if  $r=7$, and  it has
one base point, determining a rational map to $\PP^1$ if $r=8$. 
Smooth rational 
curves $E \subset X_r$ such that  $(E,E) =-1$ and $(E, -K) =1$ are called 
{\em exceptional curves}. 

\begin{theo} \cite{Ma}
The exceptional curves on $X_r$ are the following: 

\begin{enumerate}
\item blown-up points $p_1, \ldots, p_r$;
\item lines through pairs of points $p_i, p_j$;
\item conics through 5 points from  $\{p_1, \ldots, p_r\}$($r\geq5$);
\item cubics, containing 7  points and 1 of them double ($r\geq7$);
\item quartics, containing 8 points and 3 of them double ($r=8$);
\item quintics, containing 8 of point and  6 of them double ($r=8$);
\item sextics, containing 8 of those points, 7 of them double and 1 triple
($r=8$). 
\end{enumerate}
The number $N_r$ of exceptional curves on $X_r$ is given by the 
following table:
\medskip 
\begin{center} 
{\rm 
\begin{tabular}{|c|c|c|c|c|c|c|} \hline
$r$ & 3 & 4 & 5 & 6 & 7 & 8 \\ \hline
$N_r$ & 6 & 10 & 16 & 27 & 56 & 240 \\
\hline 
\end{tabular} }
\end{center} 
\end{theo} 
\medskip

The root system $R_r \subset \Pic(X_r)$ is defined  as 
\[ R_r:= \{\alpha  \in \Pic(X_r)\; : \; (\alpha,\alpha ) = -2, \; 
(\alpha, -K) = 0 \}. \]
It is easy to show that $ R_r$ is exactly the set of all classes 
$\alpha = [E_i] - [E_j]$ where 
$E_i$ and $E_j$ are two exceptional curves on $X_r$ 
such that $E_i \cap E_j = \emptyset.$  

The corresponding Weyl group $W_r$ is generated by the reflections 
 $\sigma\;  : \; x \mapsto x+(x,\alpha)\alpha$ for $\alpha \in R_r$. 
There are so called {\em simple roots} 
$\alpha_1, \ldots, \alpha_r$ such that the corresponding reflexions 
$\sigma_1, \ldots,
\sigma_r$ form a minimal generating subset of $W_r$.  
The set of simple roots can be chosen as 
\[ \alpha_1 = l_1 - l_2, \alpha_2 = l_2 - l_3, 
\alpha_3= l_0 -l_1 - l_2 - l_3, \]
\[ \alpha_i = l_{i-1} - l_i, \; \; i \geq 4. \]

The blow up morphism $X_r \to X_{r-1}$  determines an isometric 
embedding of the
Picard lattices $\Pic(X_{r-1}) \hookrightarrow \Pic(X_r)$. 
This  induces  the embeddings 
for root systems, simple roots and Weyl groups
$W_r$. For $r\geq 3$, the Dynkin diagram of $R_r$ can be considered as the 
subgraph on the vertices $\alpha_i$ $( i \leq r)$ of the following graph: 
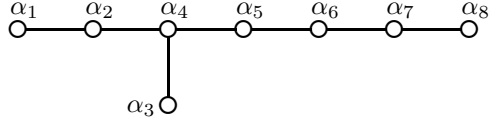
\begin{center}
\small\tt
\unitlength1cm
\begin{picture}(9,2)
\thicklines
\multiput(1,1.3)(1,0){7}{\circle{.2}}
\multiput(1.1,1.3)(1,0){6}{\line(1,0){.8}}
\put(3,.3){\circle{.2}}
\put(3,.4){\line(0,1){.8}}
\put(2.45,.2){$\alpha_3$}
\put(.9,1.5){$\alpha_1$}
\put(1.9,1.5){$\alpha_2$}
\put(2.9,1.5){$\alpha_4$}
\put(3.9,1.5){$\alpha_5$}
\put(4.9,1.5){$\alpha_6$}
\put(5.9,1.5){$\alpha_7$}
\put(6.9,1.5){$\alpha_8$}
\end{picture}
\end{center}
In particular, we obtain 
$R_3=A_2\times A_1, R_4=A_4,
R_5=D_5, R_6= E_6, R_7 = E_7, R_8 = E_8$. 

Denote by $\varpi_1, \ldots, \varpi_r $ the dual basis to 
the $\Z$-basis $- \alpha_1,  \ldots, -\alpha_r$. Each  $\varpi_i$ is 
the highest weight of an irreducible representation of 
$G(R_r)$ which  is called a {\em fundamental representation}. 
We shall denote by $V(\varpi)$ the representation space of $G(R_r)$
with the  highest weight $\varpi$.  

\begin{dfn} 
{\rm 
A dominant weight $\varpi$ is called {\em minuscule} if all weights of  
$V(\varpi)$ are nonzero and the 
$W_r$-orbit of the highest weight vector  is a 
$\Bbbk$-basis of
$V(\varpi)$ \cite{G/PI}). 
A dominant weight $\varpi$ is called {\em quasiminuscule} 
\cite{G/PIII}, if all nonzero weights  of $V(\varpi)$ have multiplicity $1$ 
and form an $W_r$-orbit of $\varpi$ (the zero weight of $V(\varpi)$ 
may have some positive multiplicity). }
\label{weights} 
\end{dfn} 

One can  see from the explicit description of the root systems $R_r$ 
that $\varpi_r$ is minuscule for $3 \leq r \leq 7$, and  $\varpi_8$ 
is quasiminuscule.

The dimension $d_r$ of of the irreducible 
representation $V(\varpi_r)$ of $G(R_r)$ is given by the 
following table:
\medskip 
\begin{center} 
{\rm 
\begin{tabular}{|c|c|c|c|c|c|} \hline
$r$  & 4 & 5 & 6 & 7 & 8 \\ \hline
$d_r$  & 10 & 16 & 27 & 56 & 248 \\
\hline 
\end{tabular} }
\end{center} 
\bigskip

We will need the following statement: 

\begin{prop}\label{BPF}
Let $D$ be a  divisor on a Del Pezzo surface $X_r$ 
$(2 \leq r \leq 8)$ such 
that $(D,E)\geq 0$ for every exceptional curve $E \subset X_r$.  Then 
the following statements hold: \\
(i) the linear system $|D|$ has no base points on any 
exceptional curve  $E \subset X_r$; \\
(ii) if  $r\leq 7$, then  the linear system $|D|$ has no base points on 
$X_r$ at all. \\
\end{prop}

\noindent
{\em Proof.} Induction on $r.$ If $r=2$, then there exists exactly $3$ 
exceptional curves $E_0, E_1, E_2$, whose classes in the standard basis 
are $l_0-l_1-l_2, l_1, l_2$. Moreover $[E_0]$, $[E_1]$ and $[E_3]$ 
form a basis of the Picard lattice $Pic(X_2)$. The dual basis
w.r.t. the intersection form is $l_0, l_0-l_1, l_0-l_2$. Therefore the 
above conditions on $D$ imply that 
$$[D] = n_0 l_0 + n_1(l_0 -l_1) + n_2 ( l_0-l_2), \; n_0, n_1, n_2 \in 
\Z_{\geq 0} $$ 
So it is sufficient to check that the linear systems with the classes
 $l_0, l_0-l_1, l_0-l_2$ have no base points. The latter immediately follows
from the fact that the first system defines the  birational morphism 
$X_2 \to \P^2$ contracting  $E_1$ and $E_2$, the second and third linear 
systems define conic bundle fibrations over $\P^1$.

For $r>2$,  we consider a second  induction on  ${\rm deg}\, D = 
(D, -K)$. 

If there is an
exceptional curve $E \subset X_r$ with $(D,E)=0,$ then 
the invertible sheaf $\OO(D)$ is the inverse image of an invertible sheaf
$\OO(D')$ on the Del Pezzo surface $X_{r-1}$ obtained by 
the contraction of $E.$ Since the pull-back of any  exceptional curve 
on $X_{r-1}$ under the birational morphism $\pi_E\; : \; 
X_r \to X_{r-1}$ is 
again an exceptional curve on $X_r$, we obtain that $D'$ satisfy all 
conditions of the proposition on $X_{r-1}$. By the induction 
assumption   ($r-1 
\leq 7$), $|D'|$ has no base points on  $X_{r-1}$. Therefore, 
$|D| = |\pi_E^*D'|$ 
has no base points on $X_r$.  

If there is no exceptional curve  $E \subset X_r$ with $(D,E)=0,$ then 
we denote by $m$ the minimal intersection number $(D,E)$ where $E$ runs 
over all exceptional curves. Since we have $(E, -K) = 1$ for 
all exceptional curves, the divisor $D':= D + mK$ has nonnegative 
intersections with all exceptional curves and there exists an exceptional
curve $E \subset X_r$ with $(D',E)=0$. Since ${\rm deg}\,D' = (D', -K) = 
(D,-K) -
 m(K,K) < (D,-K) = {\rm deg}\,D$, by the induction 
assumption, we obtain that $|D'|$ is base point free. 
If  $r\leq 7$, then the anticanonical linear system $|-K|$ has no 
base points. Therefore, $|D| = |D' + m(-K)|$ is also base point free. 
In the case $r=8$,  $|-K|$ does have a base point $p \in X_8$. However, 
$p$ cannot lie on an exceptional curve $E$, because the short
exact sequence 
\[ 0 \to  H^0(X_8, \OO(-K - E)) \to H^0(X_8, \OO(-K)) \to 
H^0(E, \OO(-K)|_E) \to 0 \]
induces an isomorphism $H^0(X_8, \OO(-K)) \cong  H^0(E, \OO(-K)|_E)$ (since 
$\deg(-K - E) =0$ and $H^0(X_8, \OO(-K-E)) = 0$).
\qed

\section{Generators of $\Cox(X_r)$}

Let $\{ E_1, \ldots, E_{N_r} \}$ be the set of all exceptional curves 
on a Del Pezzo surface $X_r$. We choose a $\Bbbk$-rational point 
$p \in U:= X_r \setminus (\bigcup_{i =1}^{N_r} E_i)$ and denote the 
ring $\Cox(X_r, U, p)$ (see \ref{cox-ring}) simply by $\Cox(X_r)$.

The ring 
$$\Cox(X_r) = \bigoplus_{[ D] \in M_{\rm eff}(X_r)} \Cox(X_r)^{[D]}
$$ is graded by the semigroup $M_{\rm eff}(X_r) \subset 
\Pic(X_r)$ 
of classes $[D]$ of effective divisors $D$ on $X_r$. 
There is a
coarser grading on $\Cox(X_r)$ given by 
$$\Cox(X_r)^d :=  
\bigoplus_{{\rm deg}\,[ D] =d}\Cox(X_r)^{[D]},$$ where 
${\rm deg}\,[ D] :=(D,-K)$.

\begin{prop}\label{toric}
The graded ring $\Cox(X_3)$ is isomorphic to  a
polynomial ring in 6 variables $\Bbbk[x_1,\dots,x_6]$, 
where $x_i$ are sections defining all $6$ 
exceptional curves on $X_3$.
\end{prop}

\noindent
{\em Proof.} The Del Pezzo  surface $X_3$ is a toric variety which can be 
described as  the blow-up of $3$ torus invariant points (1:0:0), (0:1:0) 
and (0:0:1) in $\P^2$. So we can apply a general result of Cox 
on toric varieties \cite{Cox} (see also \ref{cox-tor}). 
 \qed

\begin{theo} \label{Gen}
For $3\leq r \leq 8$,  the ring $\Cox(X_d)$ is generated by elements of 
degree 1. If $r\leq 7$, then the generators of  $\Cox(X_d)$ are 
global sections of invertible sheaves 
defining the  exceptional curves.  If  $r=8,$ then 
we should add to the above set of generators  two linearly independent  
global sections of 
the anticanonical 
sheaf on $X_8$. 
\end{theo}

\noindent
{\em Proof.} 
Induction on $r$. The case $r=3$ is settled by  the previous proposition. 

For $r> 3$ we choose an effective divisor $D$ on $X_r$. We call 
a  section $s \in H^0(X_r, \OO(D))$ a {\em distinguished global section}
if its support is contained in the union of exceptional curves of $X_r$ 
$(r \leq 7)$, or if  its support is contained in the union of 
exceptional curves of $X_8$ 
and some anticanonical curves on $X_8$. Our purpose is to show   
that the vector space $H^0(X_r, \OO(D))$ is spanned by all 
distinguished global sections. 

This will be proved by induction on ${\rm deg}\,D:= (D, -K) > 0$. 

We consider several cases:
\begin{itemize}
\medskip

\item
If there exists  an exceptional curve $E$ such  that $(D,E)<0$, then 
$ H^0(X_r, \OO(D)|_E) = 0$ and it 
follows from the exact sequence
$$ H^1(X_r, \OO(D)|_E) \to H^0(X_r, \OO(D-E))  \to H^0(X_r, \OO(D)) \to 0 $$ 
that the multiplication
by  a  non-zero distinguished global section of $\OO(E)$ 
induces an epimorphism 
$ H^0(X_r, \OO(D-E))  \to H^0(X_r, \OO(D))$. Since 
${\rm deg}\,(D-E) = {\rm deg}\,D -1$, using the induction assumption for 
$D' = D -E$, we obtain the required statement for $D$. 
\medskip

\item
If there exists  an exceptional curve $E$ such  that $(D,E)=0,$ then 
$\OO(D)$ is the inverse
image of a sheaf $\OO(D')$ on the Del Pezzo surface $X_{r-1}$ obtained 
by the contraction of $E$. Therefore 
we have an isomorphism $ H^0(X_r, \OO(D)) \cong  H^0(X_{r-1}, \OO(D'))$ 
and, by   the induction assumption for $r-1$, 
we obtain the required statement for $D$, because  distinguished 
global sections of  $\OO(D')$ lift to  distinguished 
global sections of  $\OO(D)$. 
\medskip

\item
If $D=-K,$ (or, equivalently, if $(D,E)=1$ for every exceptional curve $E$), 
then  $\OO(D)|_E)$ is isomorphic to $\OO_E(1)$ and we have 
$ H^1(X_r, \OO(D)|_E) = 0$ together with 
the   exact sequence
\[ 0   \to H^0(X_r, \OO(D-E))  \to H^0(X_r, \OO(D)) \to  
H^0(X_r, \OO(D)|_E) \to  0, \]
where $H^0(X_r, \OO(D)|_E)$ is 
$2$-dimensional. Since ${\rm deg }\,(D-E) = {\rm deg }\,D -1 < {\rm deg }\,$, 
we can  apply  the induction assumption for $D' = D -E$. It 
remains show that there exists 
two linearly independent distinguished global sections of $ \OO(D)$ 
such that their restriction to $E$ are 
two  linearly independent global sections of $ \OO(D)|_E$. 
We describe these two distinguisched sections explicitly for 
each value of $r \in \{ 4,5,6,7,8 \}$. Without loss of generality
 we can assume that $[E] = l_1$. 

If $r =4$, then we write the anticanonical class $-K =
 3l_0 - l_1 - \cdots - l_4$ in the following 
two ways: 
\begin{eqnarray*} -K & = & (l_0 - l_1 - l_2) + 
(l_0 - l_3 - l_4) + (l_0 - l_2 - l_3) + l_2 + l_3 \\
& = &
 (l_0 - l_1 - l_3) + 
(l_0 - l_2 - l_4) + (l_0 - l_2 - l_3) + l_2 + l_3. 
\end{eqnarray*}

These two decompositions of 
$-K$ determine two distinguished global  
sections of $\OO(-K)$ with support on $5$ 
exceptional curves.
The projections of these sections  
under the morphism  $X_4 \to 
\P^2$ are shown below in Figure 1.

\begin{figure}[htbp]
\begin{center}
\small
\begin{picture}(6,2)
\multiput(2,0)(3,0){2}{\circle{.2}}
\multiput(.5,1.5)(3,0){2}{\circle{.2}}
\multiput(.5,0)(3,0){2}{\usebox{\fourpts}}
\multiput(2.3,-.3)(3,0){2}{\line(-1,1){2.1}}
\put(0,0){\line(1,0){2.5}}
\put(0,1.5){\line(1,0){2.5}}
\put(.3,-.35){$E$}
\put(3.1,-.1){$E$}
\put(2.8,.6){$=$}
\put(3.5,-.5){\line(0,1){2.5}}
\put(5,-.5){\line(0,1){2.5}}
\end{picture}
\end{center}
\caption{Two distinguished  anticanonical classes for $r=4.$}
\end{figure}
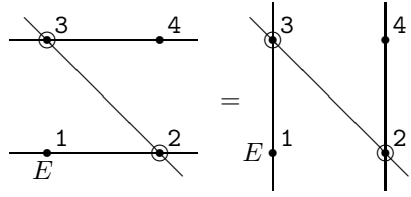

The restriction of the first section to $E$   
vanishes at the intersection point $q_1$ of $E$ with the exceptional 
curve with the class $l_0 - l_1 - l_2$. The restriction of the second section 
to $E$ vanishes at the intersection point $q_2$ of $E$ with the exceptional 
curve with the class $l_0 - l_1 - l_3$. It is clear that $q_1 \neq q_2$. So 
the distinguished  anticanonical sections are linearly independent.
\medskip

If $r=5$, then we write the anticanonical class 
as  

\begin{eqnarray*} 
-K & = &  3l_0 - l_1 - \cdots - l_5 \\ 
& = &  (l_0 - l_1 - l_2) + (l_0 -l_3 -l_4) + (l_0 - l_4  -l_5) + l_4 \\
& = &
 (l_0 - l_1 - l_5) + (l_0 -l_2 -l_3) + (l_0 - l_3  -l_4) + l_3. 
\end{eqnarray*} 
 
The corresponding  distinguished 
anticanonical sections vanish at two 
different intersection points of $E$ with 
the exceptional curves belonging to the  classes $l_0 - l_1 - l_2$ and 
$l_0 - l_1 - l_5$.  
\medskip 

If $r=6$, then we write the anticanonical class   as  
\begin{eqnarray*} -K & = & 
 3l_0 - l_1 - \cdots - l_6 \\ & = & 
 (l_0 - l_1 - l_2) + (l_0 -l_3 -l_4) + (l_0 - l_5  -l_6) \\ &=& 
 (l_0 - l_1 - l_6) + (l_0 -l_5 -l_4) + (l_0 - l_3  -l_2).
\end{eqnarray*} 

The corresponding  distinguished 
anticanonical sections vanish at two 
different intersection points of $E$ with 
the exceptional curves belonging to the classes $l_0 - l_1 - l_2$ and 
$l_0 - l_1 - l_6$.  
\medskip 

If $r=7$, then we write the anticanonical class  as  

\begin{eqnarray*} -K  & = & 3l_0 - l_1 - \cdots - l_7 \\ & = &
(2l_0 - l_1 - l_2 -l_3 -l_4 - l_5) + (l_0  -l_6 - l_7)\\ & = &
 (2l_0 - l_7 - l_6 -l_5 -l_4 - l_3) + (l_0   -l_2 - l_1) . 
\end{eqnarray*}

The corresponding  distinguished 
anticanonical sections vanish at two 
different intersection points of $E$ with 
the exceptional curves belonging to the classes 
$2l_0 - l_1 - l_2 -l_3 -l_4 - l_5$ and 
$l_0 - l_2 - l_1$. 
\medskip

If $r=8$, then ${\rm deg}\, D- E = 0$. Therefore, $H^0(X_8, \OO(D-E)) =0$ 
(see the proof of \ref{BPF})  
and we have an isomorphism 
\[ H^0(X_8, \OO(D)) \cong   
H^0(X_8, \OO(D)|_E). \]
So $H^0(X_8, \OO(D)|_E)$  is generated by the restrictions of the 
anticanonical sections and we're done.
\medskip 

\item
If $(D,E)\geq 1$ for all exceptional curves $E$ and $D \neq -K$, then 
we denote by $m$ the minimum of the numbers $(D, E)$ for all exceptional 
curves. Let $E_0$ be an exceptional curve such that $(D, E_0) = m \geq 1$. 
We define $D' = D - E_0$ and $D'':= D + mK$. 
By \ref{BPF}, $|D'|$ and $|D''|$ 
have no base points  (if $r \leq 7$). In particular, $D''$ is represented 
by an effective divisor. Since ${\rm deg}\, D'' = {\rm deg}\, D -m (K,K) < 
{\rm deg}\, D$,  $D''$ can be seen  as zero of a distinguished 
global section $s \in H^0(X_r, \OO(D+mK))$ 
whose support does not 
contain the  exceptional curve $E_0$ (if $r \leq 8$). We have the short exact
sequence  
\[ 0  \to H^0(X_r, \OO(D'))  \to H^0(X_r, \OO(D)) \to  
H^0(X_r, \OO(D)|_{E_0}) \to  0.\]
By the induction assumption, the space  $H^0(X_r, \OO(D'))$ is generated 
by distinguished global sections. It remains to show that there 
exist distinguished global sections 
of  $\OO(D)$ such that  
their restriction to $E_0$ generate the space $H^0(X_r, \OO(D)|_{E_0})$. 
Since $(-mK, E_0) = (D, E_0) = m$, the space $H^0(X_r, \OO(D)|_{E_0})$ 
is isomorphic to $H^0(X_r, \OO(-mK)|_{E_0})$. Since $(D'',E_0) =0$ the 
distinguished global section  $s \in H^0(X_r, \OO(D+mK))$ is nonzero 
at any point of $E_0$. Therefore 
the multiplication by the distinguished global section $s$  
defines a  homomorphism 
\[   H^0(X_r, \OO(-mK)) \to  H^0(X_r,\OO(D)) \]
whose restriction to  $E_0$ is an isomorphism 
$$H^0(X_r, \OO(-mK)|_{E_0}) 
\cong H^0(X_r, \OO(D)|_{E_0}.$$
Therefore, it is enough to show that restrictions of the distinguished 
global sections of $\OO(-mK)$ to $E_0$ generate the space 
$H^0(X_r, \OO(-mK)|_{E_0})$.  Our previous considerations 
have shown this for $m =1$. The general case $m \geq 1$ follows now 
immediately 
from the fact that the homomorphism $H^0(X_r, \OO(-K))\to 
H^0(E_0, \OO_{E_0}(1))$ is surjective and the space 
$H^0(E_0, \OO_{E_0}(m))$ 
is spanned by tensor products of $m$ elements from 
$H^0(E_0, \OO_{E_0}(1))$. \qed
\end{itemize}

\begin{coro} 
The semigroup $M_{\rm eff}(X_r) \subset \Pic(X_r)$ of classes 
of effective divisors on 
a Del Pezzo surfaces  $X_r$ $( 2 \leq r \leq X_r)$ is generated by 
elements of degree $1$. These elements are exactly the classes of 
exceptional curves if $r \leq 7$ and the classes of 
exceptional curves together with the anticanoncal class for $r =8$. 
\label{monoid}
\end{coro} 
 
\begin{prop}
If $D$ is an effective divisor of degree $\geq 2$ on $X_8$, then 
the vector space $H^0(X_8, \OO(D))$ is spanned by distinguished 
global sections of $\OO(D)$ with supports only on exceptional curves. 
\end{prop} 

\noindent
{\em Proof.} By \ref{Gen} and \ref{monoid}, it is sufficient 
to check the statement for $D = -2K$ and for $D = -K + E$ for 
any exceptional curve. The latter case immediately follows from 
\ref{Gen}, because  $D = -K +E$ is the pull back of 
the anticanonical sheaf on $X_7$ obtained by the contraction of $E$.
In the case $D = -2K$, we obtain  $120$ distinguished global sections
of $\OO(D)$ from $120$ pairs of exceptional curves $E_i, E_j$ such that 
$(E_i, E_j) = 3$:
\[ -2K = 6l_0 - 2l_1 - \ldots - 2l_8  = l_1 + 
(6l_0 - 3l_1 - 2l_2 \ldots - 2l_8). \]  
It is well-known (see e.g. \cite{Dem})
 that $X_8$ can be realized as a hypersurface of 
degree $6$ in the weighted projective space $\P(3,2,1,1)$. In particular, 
the linear system $|-2K|$ defines a double covering of $X_8$ over 
a singular quadratic cone ${\mathcal Q} \cong \P(2,1,1) \subset \P^3$. 
The single singular point $p \in
Q$ is the image of the base-point $b \in X_8$ of $|-K|$ on $X_8$. Let $C 
\subset {\mathcal Q}$ be the ramification locus ($C$ is a curve of degree $6$ 
in  $\P(2,1,1)$).  Then $120$ pairs of  
exceptional curves  $E_i, E_j$ on $X_8$ 
such that $[E_i] + [E_j] = 2[-K]$ one-to-one correspond to conics in  
$\P(2,1,1)$ which are $3$-tangent to the ramification curve $C$. 
Since every such conic in ${\mathcal Q}$ is uniquely determined as 
${\mathcal Q} \cap H$ for some plane $H \subset \P^3$. Therefore, the 
distinguished sections in $H^0(X_8, \OO(-2K))$ can be identified (up to 
a scalar multiple) with the above planes $H \subset \P^3$. 
It remains to show that 
all these $120$ planes $H$ cannot pass through the some common point 
$x \in \P^3$ for a
generic choice of the sextic $C \subset \P(2,1,1)$. 
The later can be checked by standard dimension arguments. 
\qed

\begin{rem}
\label{weigh}
Since $H^0(X_r, \OO(E))$ is $1$-dimensional for 
each exceptional curve $E \subset X_r$, we can choose a nonzero section 
$x_E \in  H^0(X_r, \OO(E))$ which is determined  up to multiplication
by a nonzero scalar. Therefore the affine algebraic variety 
$\A(X_r):= {\Spec}\,\Cox(X_r)$ is embedded into 
the affine space $\A^{N_r}$ on which the maximal 
torus $T_r \subset G(R_r)$ acts in a canonical way 
such that the space  $\A^{N_r}$ can 
identified with the  representation space $V(\varpi_r)$ 
of the  algebraic group $G(R_r)$  (if  $ r \leq 7$). 
In the case $r=8$, all $240$ exceptional curves on $X_8$ can be similarly 
identified with all non-zero weights
of the adjoint representation of $G(E_8)$ in $V(\varpi_8)$.
The space $V(\varpi_8)$  contains   
the weight-$0$  subspace of dimension $8$, but the ring $\Cox(X_r)$ has only 
2-dimensional space of anticanonical sections. Thus, we cannot 
identify the  degree-$1$ homogeneous component of   
$\Cox(X_8)$ with the  representation space  $V(\varpi_8)$ of $G(E_8).$
\end{rem}

Since the kernel of the surjective homomorphism
\[ \deg \, : \, Pic(X_r) \to \Z, \]
can be identified with the character group ${\mathfrak X}(T_r)$ of a maximal 
torus $T_r \subset G(R_r)$ and the torus  $T_r$ acts on the homogeneous
space $G(R_r)/P(R_{r-1})$ embedded into the projective space 
$\P V(\varpi_r)$ we obtain a natural $\Pic(X_r)$-grading of 
the homogeneous coordinate ring 
 of the projective variety 
$G(R_r)/P(R_{r-1})$. 

\begin{theo} 
Let  $\lambda$ be an element in $\Pic(X_r)$.  
The weight-$\lambda$ subspace in the homogeneous coordinate ring 
 of the projective variety 
$G(R_r)/P(R_{r-1})$ is nonzero if and only if $\lambda$ is 
represented by  an effective divisor on $X_r$ (i.e., $\lambda \in 
M_{\rm eff}(X_r)$). 
\label{grading}
\end{theo} 

\noindent
{\em Proof.} It is known that  the projective variety 
$G(R_r)/P(R_{r-1})$ is arithmetically normal and Cohen-Macaulay \cite{DL,
G/PV}. In particular, the  homogeneous coordinate ring of $G(R_r)/P(R_{r-1})$
is generated 
by elements of degree $1$.  Therefore, the  weight-$\lambda$ subspace
in the cooordinate ring is nonzero if and only if $\lambda$ is 
a nonnegative integral linear combination of $\Pic(X_r)$-weights having 
positive 
multiplicity in   $V(\varpi_r)$. 
By \ref{monoid} and \ref{weigh}, the latter is 
equivalent to  $\lambda \in 
M_{\rm eff}(X_r)$.

\section{Quadratic  relations in $\Cox(X_r)$}

Let us denote  $\P(X_r) := \Proj\,\Cox(X_r)$. 
If  $4 \leq r \leq  7$, then the projective variety $\P(X_r)$ is canonically 
embedded into the  projective space $\P^{N_r -1}$ ($N_r$ is 
the number of exceptional curves on $X_r$). The affine variety
$\A(X_r) \subset \A^{N_r}$ is the affine cone over $\P(X_r)$.

\begin{prop} 
The ring $\Cox(X_4)$ is isomorphic to the subring of all $3\times 3$-minors
of a generic $3 \times 5$-matrix. In particular, 
the projective variety $\P(X_4) \subset \P^9$ is isomorphic 
to the Pl{\"u}cker embedding of the Grassmannian 
$Gr(3,5)$. 
\label{G35}
\end{prop}

\noindent
{\em Proof.} In order to describe the multiplication in 
$\Cox(X_4),$ one needs to choose a basis in
$\Pic(X_4).$ 

Let
$x:y:z$ be the homogeneous coordinates on $\PP^2$. 
We choose the basis $l_0,\dots,l_4,$ as in Section 2, i.e., 
$l_0$ is  the preimage of the line $z =0$ at infinity, 
$l_1, l_2, l_3, l_4$ are classes of 
the exceptional fibers over $4$ points $p_1, \ldots, p_4 \in 
\PP^2$.
We identify the representatives of each  class in $\Pic(X_4)$ with
the subsheaves $\OO(\sum_{i =0}^4 m_il_i)$ of the 
constant sheaf $\KK(X_4)$ of rational functions on $X_4.$ Then
the ring multiplication in $\Cox(X_4)$ is just the multiplication of the
corresponding rational functions in  $\KK(X_4)$.

Let $(x_i:y_i:z_i)$  be the coordinates of the blown-up point 
$p_i \in \P^2$ ($i=1,
\dots,4$). Consider the $3\times 5$-matrix
$$M=\left(
\begin{array}{ccccc}
x_1&x_2&x_3&x_4&x/z\\
y_1&y_2&y_3&y_4&y/z\\
z_1&z_2&z_3&z_4&1
\end{array}
\right).
$$

For any $3$-element subset $I= \{ i,j, k\}  \subset 
\{ 1, \ldots, 5 \}$, we denote 
by  $M_{I}$  the maximal minor of $M$
consisting of the columns with numbers in $I$ taken in the natural order.

We choose the rational functions in  $\KK(X_4)$ representing the generators  
$x_E$ of $\Cox(X_4)$ as follows: 

\[ x_{l_1}:=M_{\{2,3,4\}}, \; x_{l_1}:=M_{\{1,3,4\}}, \; x_{l_3}:=M_{\{1,2,4\}}, \;x_{l_1}:=M_{\{1,2,3\}},  \]
\[ x_{l_0 - l_i - l_j}:= M_{\{i,j,5\}}, \;\;  1\leq i < j \leq 4. \]
All these functions are non-zero because the points $p_1, \ldots, p_4$ 
are in general position.

It is known that the generators of the homogeneous coordinate ring of $G(3,
5)$ are naturally identified with the maximal minors of a generic $3\times
5$-matrix. 
Consider the homomorphism $\varphi$ of the homogeneous
coordinate ring of $G(3,5)$ to $\Cox(X_4),$
which  sends these generic minors into the corresponding minors of the
matrix $M$ above. Since $\Cox(X_4)$ is generated by 
$\{x_E\}$, this homomorphism 
is surjective.
By \ref{grading}, $\varphi$  respects the $\Pic(X_r)$-grading (in particular, 
$\varphi$ respects the $\ZZ_{\geq 0}$-grading as well). The surjectivity 
of $\varphi$  induces a closed embedding of $\P(X_4)$ into $G(3,5).$ 
Since both varieties
are irreducible of dimension 6 (see \ref{dimension}), we obtain an isomorphism of $\P(X_4)$ and $G(3,5)$ as subvarieties of $\P^9$. Therefore 
$\varphi$ is an isomorphism of the homogeneous coordinate ring of $G(3,5)$ 
and  $\Cox(X_4)$. In particular, $\Cox(X_4)$ is defined by 
$5$ quadratic Pl{\"u}cker relations. One of these relations is 
\[ M_{\{1,2,5\}} M_{\{3,4,5\}}  -  M_{\{1,3,5\}} M_{\{2,4,5\}} +
 M_{\{1,4,5\}} M_{\{2,3,5\}} =0. \]
\qed
\bigskip

The article \cite{G/PI} describes a $\Bbbk$-basis for the homogeneous
coordinate ring of $G/P$ in the case, when $P$ is a maximal parabolic
subgroup containing a Borel subgroup $B$ such that the fundamental
weight $\varpi$ corresponding to $P$ is
minuscule (see \ref{weights}). It also shows that this ring is always 
defined by quadratic relations.

A way to write explicitly the quadratic 
relations for the orbit of the highest weight
vector for any representation of a semisimple Lie group is given in
\cite{Li}. A more geometric approach to these quadratic equations
is contained in the proof of Theorem 1.1 in \cite{L-T}:

\begin{prop}\label{S2R}
The orbit $G/P_{\varpi}$ of the highest weight vector in the projective space 
$\PP V(\varpi)$ is the
intersection of the second Veronese embedding of $\PP V(\varpi)$
with the subrepresentation $ V(2\varpi)$ of the symmetric square
$S^2 V(\varpi)$. Moreover, these quadratic relations generate the 
ideal of  $G/P_{\varpi} \subset\PP V(\varpi)$.  
\end{prop}

We expect that the following general statement is true:

\begin{conj} 
The ideal of relations between the degree-$1$ generators of 
$\Cox(X_r)$ is generated by quadrics for  $4 \leq r \leq 8$. 
\end{conj}

 For any exceptional curve
$E \subset X$, we consider the  open chart $U_E \subset \P^{N_r -1}$ defined 
by the condition $x_E \neq 0$. Thus, we obtain an open covering 
of  $\P(X_r)$  by $N_r$ affine subsets  $U_E \cap \P(X_r)$.

\begin{prop}\label{cones}
Let $X_{r-1}$ the Del Pezzo surface obtained by the contraction 
of $E$ on $X_r$. Then there exist a natural isomorphism
\[ U_E \cap \P(X_r) \cong \A(X_{r-1}). \]
\end{prop} 

\noindent
{\em Proof.} 
Let $\pi \; : \; X_r \to X_{r-1}$ be the contraction of $E$. Then we obtain 
the ring homomorphism $\pi^* \; : \;\Cox(X_{r-1}) \to\Cox(X_{r})$.  
We shall show that the localization  $\Cox(X_{r})_{x_{E}}$ of the ring 
$\Cox(X_{r})$ by the element $x_E$ can be identified with 
the Laurent polynomial 
extension of $\pi^*\Cox(X_{r-1})$ by  $x_E$, i.e. there exist a ring 
isomorphism
\[\Cox(X_{r})_{x_{E}} \cong  \pi^*\Cox(X_{r-1})[ x_E, x_E^{-1} ] . \]
For simplicity, we assume that $[E] = l_r$ and $\{ l_0, \ldots, 
l_{r-1} \}$ is  the pull-back of the standard basis in 
$\Pic(X_{r-1})$. 
We remark that 
any divisor class $$[D] = m_rl_r  + \sum_{i =0}^{r-1} m_i l_i  
\in \Pic(X_r)$$ is  uniquely represented
as sum $m_r [E] + [D']$ where $[D']  =  \sum_{i =0}^{r-1} m_i l_i  
\in  \pi^*(\Pic(X_{r-1}))$
Using 
$\pi^*$, we identify two fields of rational functions 
${\mathcal K}(X_{r-1})$ and ${\mathcal K}(X_{r})$. This identification 
allows us to consider 
the vector space 
\[  H^0( X_r, \OO( m_rl_r  + \sum_{i =0}^{r-1} m_i l_i)) \]
 as a subspace of 
\[  H^0( X_{r-1}, \OO(\sum_{i =0}^{r-1} m_i l_i)) x_E^{m_r} \subset 
\pi^*\Cox(X_{r-1})[ x_E, x_E^{-1}]. \]
For fixed integers $m_0, m_1, \ldots , m_{r-1}$, the 
embedding of the vector spaces 
\[  H^0( X_r, \OO( m_rl_r  + \sum_{i =0}^{r-1} m_i l_i))  \hookrightarrow 
 H^0( X_{r-1}, \OO(\sum_{i =0}^{r-1} m_i l_i)) x_E^{m_r} \]
is an isomorphism for sufficiently large $m_r$. 
Moreover, this embedding of vector spaces respects the multiplications 
in  $\Cox(X_{r})$ and in $\pi^*\Cox(X_{r-1})[ x_E, x_E^{-1} ]$.  
Thus, we obtain an embedding of rings
\[  \Cox(X_{r}) \hookrightarrow  \pi^*\Cox(X_{r-1})[ x_E, x_E^{-1} ]. \]
On the other hand, it is clear that $\pi^*\Cox(X_{r-1})[ x_E, x_E^{-1}]$ 
is a subring 
of the localization $\Cox(X_{r})_{x_{E}}$. Thus, we get an isomorphism
\[\Cox(X_{r})_{x_{E}} \cong  \pi^*\Cox(X_{r-1})[ x_E, x_E^{-1} ] . \]

Now we remark that 
the coordinate ring of the affine variety  $U_E \cap \P(X_r)$ 
is degree-$0$ component of  $\Cox(X_{r})_{x_{E}}$. By the above isomorphism, 
this component is isomorphic to $\Cox(X_{r-1})$.   
\qed

\begin{coro} 
The singular locus of the algebraic varieties  $\P(X_r)$ and  $\A(X_r)$ 
has  codimension $7$. 
\end{coro} 

\noindent
{\em Proof.} Since $\A(X_3) \cong \A^6$, we obtain that  $\P(X_4)$ is a smooth
variety covered by $10$ affine charts which are isomorphic to  $\A^6$. 
Using the isomorphism $\P(X_4) \cong G(3,5)$ (see \ref{G35} ), we obtain 
that $\A(X_4)$  has an isolated singularity at $0$. Therefore, the singular 
locus of $\P(X_5)$ consists of $16$ isolated points. The singular 
locus of $\P(X_6)$ is $1$-dimensional and the singular locus of 
$\P(X_7)$ is $2$-dimensional etc. \qed

\begin{dfn}\label{rul}
{\rm A divisor class $[D]$ is called \emph{a ruling} if it can be written 
as  a sum of two classes of exceptional curves 
$[E_i] + [E_j]$ such that $(E_i, E_j) =1$, or , equivalently, if $D$ satisfies 
the conditions $(D,D)=0,\ (D,-K)=2.$  The 
invertible sheaf corresponding to a ruling  determines a conic 
bundle morphism $X_r\to \PP^1.$}
\end{dfn}

\begin{rem}
Lemma~5.3 of~\cite{F-M} says that the Weyl group acts transitively on
rulings. 
\end{rem} 

Each ruling  $[D]$ can be represented by  $r-1$ different ways 
as   a sum of two classes of exceptional curves corresponding to degenerate
fibers of the conic bundle  $X_r\to \PP^1.$. Thus, we obtain 
$r-1$ distinguished sections in the $2$-dimensional space 
$H^0(X_r, \OO(D))$. If  $r \geq 4$, then  for each ruling $[D]$, 
 we obtain in this way $r-3$ 
linearly independent quadratic relations between generators of  
$\Cox(X_r)$. 

\begin{rem}
We note that  $Pic(X_4)$ has exactly $5$ rulings.  Each such a ruling 
defines a Pl{\''u}cker quadric (see the proof of \ref{G35}).
\label{rul4}
\end{rem} 
  
We cannot expect in general that all quadratic relation among generators 
are coming from rulings. However, the following statement is true:

\begin{theo}
For $4 \leq r \leq 6$, the ring $\Cox(X_r)$ is defined by
the radical of the ideal generated by the
quadratic relations corresponding to rulings.
\end{theo}

\noindent 
{\em  Proof.}
Let $Z_r \subset \A^{N_r}$ is the affine subvariety defined by the 
quadratic relations coming from rulings. We want to show that 
$Z_r = \A(X_r)$ $(4  \leq r \leq 6$). 

For $r = 4$, the statement
follows from \ref{rul4}. 

Obviously, the zero $0 \in \A^{N_r}$ is 
common point of $Z_r$ and  $\A(X_r)$ for all $r$. 
Consider the affine open coverings of  $Z_r \setminus \{0\}$ and  $\A(X_r)
\setminus \{0\} $ defined by 
affine open subsets $x_E \neq 0$, where $E$ runs over all exceptional
curves of $X_r$. 
Using the induction on $r$ and Proposition \ref{cones}, 
we want to show that $Z_r \cap U_E = \A(X_r) \cap U_E$ for each 
exceptional curve. For this purpose, it is important 
to remark that 
 the affine coordinate
ring of $Z_r \cap U_E $ is generated
by all elements $x_F/x_E$ such that $(E,F)=0.$ 
For $r = 5,6$, the last property follows from the fact that if $(E,E') >0$ for 
two exceptional curves $E,E'$ on $X_r$, then $(E,E') =1$, i.e., $[E] + [E']$ is a 
ruling and there exists a ruling quadratic relation 
\[ x_E x_{E'} = \sum_{i} a_i X_{E_i}X_{E_i'}, \]
where all exceptional curves $E_i, E_i'$ do not intersect $E$.  
The last property shows that 
\[ Z_r  \cap U_E \cong Z_{r-1} \times (\A^1 \setminus \{0\}). \] 
It follows from the proof of  \ref{cones} that 
\[  \A(X_r)  \cap U_E \cong \A(X_{r-1}) \times (\A^1 \setminus \{0\}). \]
By induction, we have the equality $ Z_{r-1} = \A(X_{r-1})$. This implies 
the equality   $Z_r \cap U_E = \A(X_r) \cap U_E$ for each 
exceptional curve.  Thus,   $Z_r = \A(X_r)$.
\qed

\end{document}